\documentclass[11pt]{article}
\usepackage{amssymb,amsmath,amsthm}
\usepackage[english]{babel}
\usepackage{t1enc}
\usepackage[latin2]{inputenc}
\usepackage{epsfig}

\newtheorem{thm}{Theorem}%[section]

\newtheorem{lem}[thm]{Lemma}

\newtheorem{conj}[thm]{Conjecture}

\usepackage{indentfirst}
\def\Z{\mathbb Z}

\begin{document}

\title{Almost optimal pairing strategy for Tic-Tac-Toe\\ with numerous directions}
\author{Padmini Mukkamala and D\"om\"ot\"or P\'alv\"olgyi
%\footnote{ELTE, Budapest and EPFL, Lausanne. %\newline  Ecole Polytechnique Fédérale de Lausanne
%Supported by OTKA NK 67867.}
}
%\date{}
\maketitle

\begin{abstract}
We show that there is an $m=2n+o(n)$, such that, in the Maker-Breaker game played on $\Z^d$ where Maker needs to put at least $m$ of his marks consecutively in one of $n$ given winning directions, Breaker can force a draw using a pairing strategy.
This improves the result of Kruczek and Sundberg \cite{KS1} who showed that such a pairing strategy exits if $m\ge 3n$. %Maker needs to put $3n$ of his marks consecutively.
A simple argument shows that $m$ has to be at least $2n+1$ if Breaker is only allowed to use a pairing strategy, thus the main term of our bound is optimal.
%Moreover, if $2n+1$ is a prime and it works in general if $m$ is at least as big as the smallest prime bigger than $2n$.???
\end{abstract}

\medskip

\section{Introduction}
A central topic of combinatorial game theory is the study of positional games, the interested reader can find the state of the art methods in Beck's Tic-Tac-Toe book \cite{B}. In general, positional games are played between two players on a {\em board}, the fields of which they alternatingly occupy with their marks and whoever first fills a {\em winning set} completely with her/his marks wins the game.
Thus a positional game can be played on any hypergraph, but in this paper, we only consider %countable boards where
{\em semi-infinite} games
where all winning sets are finite.
If after countably many steps none of them occupied a winning set, we say that the game ended in a draw.
It is easy to see that we can suppose that the next move of the players depends only on the actual position of the board and is deterministic.\footnote{This is not the case for infinite games and even in semi-infinite games it can happen that the first player can always win the game but there is no $N$ such that the game could be won in $N$ moves. For interesting examples, we refer the reader to the antique papers \cite{ACN, BC, CN}.}
We say that a player has a {\em  winning strategy} if no matter how the other player plays, she/he always wins.
We also say that a player has a {\em drawing strategy} if no matter how the other player plays, she/he can always achieve a draw (or win).
A folklore strategy stealing argument shows that the second player (who puts {\em his} first mark after the first player puts {\em her} first mark, as ladies go first) cannot have a winning strategy, so the best that he can hope for is a draw.
Given any semi-infinite game, either the first player has a winning strategy, or the second player has a drawing strategy. %, but it might be undecidable to determine which.cite???
We say that the second player can achieve a {\em pairing strategy draw} if there is a matching among the fields of the board such that every winning set contains at least one pair. It is easy to see that the second player can now force a draw by putting his mark always on the field which is matched to the field occupied by the first player in the previous step (or anywhere, if the field in unmatched).
Note that in a relaxation of the game for the first player, by allowing her to win if she occupies a winning set (not necessarily first), the pairing strategy still lets the second player to force a draw. 
%Note that in this case the second player forfeits his option to force the first player to prevent her sometimes from completing a winning set, he can even prevent the first player from completing a winning set if she does not have to care about whether he completes one before her.
Such drawing strategies are called {\em strong draws}.
Since in these games only the first player is trying to complete a winning set and the second is only trying to prevent her from doing so, in these games, the first player is called {\em Maker}, the second {\em Breaker}, and the game is called a {\em Maker-Breaker} game.

This paper is about a generalization of the Five-in-a-Row game\footnote{Aka Go-Muku and Am\H oba.} which is the more serious version of the classic Tic-Tac-Toe game. This game is played on the $d$-dimensinal integer grid, $\Z^d$, and the winning sets consist of $m$ consecutive gridpoints in $n$ previously given directions. For example, in the Five-in-a-Row game $d=2$, $m=5$ and $n=4$, the winning directions are the vertical, the horizontal and the two diagonals.
Note that we only assume that the greatest common divisor of the coordinates of each direction is $1$, so a direction can be arbitrarily long, eg.\ $(5,0,24601)$.
The question is, for what values of $m$ can we guarantee that the second player has a drawing strategy?
It was shown by Hales and Jewett \cite{B}, that for the four above given directions and $m=9$ the second player can achieve a pairing strategy draw. In the general version, a somewhat weaker result was shown by Kruczek and Sundberg \cite{KS1}, who showed that the second player has a pairing strategy if $m\ge 3n$ for any $d$. 
They conjectured that there is always a pairing strategy for $m\ge 2n+1$, generalizing the result of Hales and Jewett.\footnote{It is not hard to show that if $m=2n$, then such a strategy might not exist, we show why in Section 3.}

\begin{conj}\cite{KS1} If $m=2n+1$, then in the Maker-Breaker game played on $\Z^d$, where Maker needs to put at least $m$ of his marks consecutively in one of $n$ given winning directions, Breaker can force a draw using a pairing strategy.
\end{conj}

Our main result almost solves their conjecture.

\begin{thm}\label{one} There is an $m=2n+o(n)$ such that in the Maker-Breaker game played on $\Z^d$, where Maker needs to put at least $m$ of his marks consecutively in one of $n$ given winning directions, Breaker can force a draw using a pairing strategy.
\end{thm}

In fact we prove the following theorem, which is clearly stronger because of the classical result \cite{H} showing that there is a prime between $n$ and $n+o(n)$.

\begin{thm}\label{two} If $p=m-1\ge 2n+1$ is a prime, then in the Maker-Breaker game played on $\Z^d$, where Maker needs to put at least $m$ of his marks consecutively in one of $n$ given winning directions, Breaker can force a draw using a pairing strategy.
\end{thm}

The proof of the theorem is by reduction to a game played on $\Z$ and then using the following recent number theoretic result of Preissmann and Mischler. Later it was discovered by Kohen and Sadofschi \cite{KS} that there is a short proof using the Combinatorial Nullstellansatz \cite{A}.

\begin{lem}\label{prime}\cite{PM} Given $d_1,\ldots,d_n$ and $p\ge 2n+1$ prime, we can select $2n$ numbers, $x_1,\ldots,x_n,y_1,\ldots,y_n$ all different modulo $p$ such that $x_i+d_i\equiv y_i \mod p$.
\end{lem}

We prove our theorem in the next section and end the paper with some additional remarks.

\section{Proof of Theorem \ref{two}}
%First we sketch the main ideas of the proof, then we give a more detailed version.

We consider the winning directions to be the primivite vectors\footnote{A vector $(v_1,\ldots,v_d)\in \Z^d$ is primitive if $gcd(v_1,\ldots,v_d)=1$.} $\vec v_1,...,\vec v_n$. 
Using a standard compactness argument it is enough to show that there is a pairing strategy if the board is $[N]^d$, where $[N]$ stands for $\{1,\ldots,N\}$. For interested readers, the compactness argument is discussed in detail at the end of this section.
%For all futher arguments,  $\in [N]^d$. And for the game to be meaningful, one should think of $N$ as large enough to contain winning sets in all directions.

First we reduce the problem to one dimension.
Take a  vector $\vec r = (r_1,r_2,...,r_d)$ and transform each grid point $\vec v$ to $\vec v\cdot \vec r$. 
If $\vec r$ is such that $r_j>0$ and $r_{j+1} > N(r_1+\ldots+r_j)$ for all $j$, then this transformation is injective from $[N]^d$ and each winning direction is transformed to some number, $d_i = |\vec r \cdot \vec v_i|$.\footnote{It is even possible that some of these numbers are zero, we will take care of this later.} 
So we have these $n$ differences, $d_1,\ldots, d_n$, and the problem reduces to avoiding arithmetic progressions of length $m$ with these differences.
From the reduction it follows that if we have a pairing strategy for this game, we also have one for the original.

%\dom{rewrite}
Let $p$ be a prime such that $2n+1 \le p \le 2n+1+o(n)$. (A classic number theoretical result \cite{H} asserts that we can always find such a $p$). 
If we pick a vector $\vec u$ % = [u_1,u_2,...,u_d]$
uniformly at random from $[p]^d$, 
then for any primitive vector $\vec v$, % \in [N]^d$,
 $\vec u \cdot \vec v$ 
%$(1\ldots p,pN+1\ldots p(N+1),p(N^2+2N)+1\ldots p(N^2+2N+1),\ldots, p(N+1)^{n-1}-p+1\ldots p(N+1)^{n-1})$, then
will be divisible by $p$ with probability $1/p$.
Since each winning direction was a primitive vector, 
using the union bound, the probability that at least one of the $\vec u \cdot \vec v_i$ is divisible by $p$ is at most $n/p<1/2$. 
So, there is a $\vec u' = (u_1',u_2',..,u_d') \in [p]^d$ such that none of $\vec u' \cdot \vec v_i$ is divisible by $p$. 
If we now take $\vec r = (r_1,r_2,..,r_d)$ such that $r_j = u_j' + (pN)^{j-1}$, then the dot product with $\vec r$ is injective from $[N]^d$ to $\Z$ and
%$$ N(r_1 + ... + r_i) < N( (p+1) + (p+ 2pN) + ... + (p + (2pN)^{i-1})$$
%$$ = ipN + \frac{(2^ip^iN^i - 1)N}{2pN-1}$$
%$$ < ipN + \frac{2^ip^iN^{i+1}}{pN} < 2^ip^iN^i \le r_{i+1}$$
%Also, this $r$ will satisfy the property that
none of the $d_i = \vec r \cdot \vec v_i$ are divisible by $p$, since $\forall j \ r_j \equiv u_j' \mod p$.  % as the probability that any of them is is less than $n/p<1/2$.

%Denote $d_i \mod p$ by $d_i'$.

We now apply Lemma \ref{prime} for $d_1,... , d_n$ to get $2n$ distinct numbers $x_1,x_2,...x_n,y_1,y_2,..,y_n$ such that $0 \le x_i,y_i < p$ and $x_i + d_i \equiv y_i \mod p$. 
Our pairing strategy is, for every $x \equiv x_i \mod p$, $x$ is paired to $x+d_i$ and if $x \equiv y_i \mod p$, then $x$ is paired to $x-d_i$.

To see that this is a good pairing strategy, consider an arithmetic progression $a_1,..., a_m$ of $m=p+1$ numbers with difference, say, $d_i$. Then, since $p$ and $d_i$ are coprimes, of the numbers $a_1,..., a_{m-1}$, one of the numbers, say $a_j$, must be such that $a_j \equiv x_i \mod p$ and hence, $a_j,a_{j+1}$ must be paired in our pairing strategy showing both cannot be occupied by Maker.\hfill$\Box$\\

For completeness here we sketch how the compactness argument goes. We show that it is sufficient to show that a pairing strategy exists for every finite $[N]^d$ board. For this we use the following lemma (note: we use the version in \cite{D}).

\begin{lem}\label{kil}\cite{KIL} (K\"onig's Infinity Lemma) Let $V_0,V_1,..$ be an infinite sequence of disjoint non-empty finite sets, and let $G$ be a graph on their union. Assume that every vertex $v$ in a set $V_N$ with $n\ge1$ has a neighbour $f(v)$ in $V_{N-1}$. Then $G$ contains an infinite path, $v_0v_1...$ with $v_N \in V_N$ for all $N$.
\end{lem}

Given a pairing strategy for $[N_0]^d$, consider a smaller board $[N]^d$ where $N<N_0$. We can think of a pairing strategy as, essentially, a partition of $[N_0]^d$ into pairs and unpaired elements\footnote{Note that a pairing strategy does not gaurentee that every element is paired. It only states that every winning set has a pair. Hence there might be many unpaired elements in a pairing strategy.}. We can construct a good pairing strategy for the smaller board by taking the restriction of these set of pairs to $[N]^d$ and leave the elements paired outside $[N]^d$ as unpaired elements. We call this as a restriction of the pairing strategy to the new board. As long as we do not change the length of the winning sets and the prescribed directions, any winning set in the $[N]^d$ board is also a winning set in the $[N_0]^d$ board and hence must have a pair from the restriction. Hence, the Breaker can block all winning pairs and the restriction of the pairing strategy is a valid strategy for Breaker for the smaller board. 

We can now prove the following theorem,
\begin{thm}\label{compact} Given a fixed set $S,\ |S|=n$, of winning directions,  and positive integer $m$, if Breaker has a pairing strategy for all boards $[N]^d$ % with $N>N_0$
 and length of winning sets equal to $m$, then Breaker also has a pairing strategy for the $\Z^d$ board. 
\end{thm}

%Let $n_0$ be the smallest integer such that the game with $S,m$ and $[n_0]^d$ board is meaningful.???
We will apply K\"onig's Infinity Lemma to prove the theorem. Let $V_N$ be the set of all pairing strategies on the $[N]^d$ board with winning sets as defined in the theorem. We say a strategy in $V_{N-1}$ and a strategy in $V_N$ have an edge between them if the former is a restriction of the latter. It is easy to see that every vertex in $V_N$ does have an edge to its restriction in $V_{N-1}$. Hence, by the lemma, we must have an infinite path $v_0v_1...$. The union of all these pairing strategies gives a valid pairing strategy for the infinite game.

%Our pairing strategy is to pair $k$ to $k+ d_i'$ if $k\equiv x_i \mod p$.
%It is easy to see that this is a good matching as $d_i$ and $p$ are coprimes.

\section{Possible further improvements and remarks}
As we said before, if $m\le 2n$, then the second player cannot have a pairing strategy draw.
This can be seen as follows. On one hand, in any pairing strategy, from any $m$ consecutive points in a winning direction, there must be at least two points paired to each other in this direction.
On the other hand, there must be a winning direction in which at most $1/n$ of all fields are matched to another in this direction.
If we pick a random set of size $m-1$ in this direction, then it will contain at most $(m-1)/n< 2$ fields matched in this direction.
Thus, there is a set of size $m-1$ that contains only one such field.
Its matching field can now be avoided by extending this set to one way or the other, thereby giving us a winning set with no matched pair.

If $n=1$ or $2$, then a not too deep case analysis shows that the first player has a winning strategy if $m=2n$, even in the normal game, where the second player also wins if he occupies a winning set.
Moreover, the second player has a pairing strategy for $m=2n+1$ if $n=1$ or $2$, thus, in this case, the conjecture is tight.
However, for higher values, it seems that Breaker can always do better than just playing a pairing strategy, so we should not expect this strategy the best to achieve a draw. Quite tight bounds have been proved for Maker-Breaker games with {\em potential} based arguments, for the latest in generalization of Tic-Tac-Toe games, see \cite{KS2}.
Despite this, from a combinatorial point of view it still remains an interesting question to determine the best pairing strategy. Unfortunately our proof can only give $2n+2$ (if $2n+1$ is a prime) which is still one bigger than the conjecture.

One could hope that maybe we could achieve a better bound using a stronger result than Lemma \ref{prime} (see for example the conjecture of Roland Bacher in their paper, we would like to thank him for directing us to it \cite{Ba}), however, already for $n=3$, our method cannot work. 
Consider the three directions $(1,0),(0,1),(1,1)$. Optimally, we would hope to map them to three numbers, $d_1,d_2,d_3$, all coprime to $6$, such that we can find $x_1,x_2,x_3,y_1,y_2,y_3$ all different modulo $6$ such that $x_i+d_i\equiv y_i \mod 6$. 
But this is impossible since $d_3=d_1+d_2$, so we cannot even fulfill the coprimeness, but even if we forget about that condition, it would still be impossible to find a triple satisfying $d_3=d_1+d_2$. 
Consider a pairing strategy where the pair of %$v$ is determined by $f(v \cdot r)$, or in other words, there is a vector $r$ such that the pair of
any grid point $\vec v$, depends only on $v \cdot r$, then the above argument shows that such a pairing strategy does not exist for the three vectors $(1,0),(0,1),(1,1)$.
%This shows that there is no pairing strategy for these three vectors satisfying the additional condition that the tile of each gridpoint $v$ is only a function of $r\cdot v$ for some $r$. 
However, it is not hard to find a suitable periodic pairing strategy. We would like to end with an equivalent formulation of Conjecture 1.

\begin{conj}[Kruczek and Sundberg, reformulated] Suppose we are given $n$ elements, $\vec v_i$ of $\Z_{2n}^{d}$ for $i\in [n]$. Is it always possible to find a partition of $\Z_{2n}^{d}$ into $\vec x_i^j,\vec y_i^j$ for $i\in [n], j\in [2n]$ such that $\vec x_i^j+\vec v_i=\vec y_i^j$ and $\vec x_i^j-\vec x_i^{j'}$ is not a multiple of $\vec v_i$ for $j\ne j'$?
\end{conj}

Also, one can formulate a more daring conjecture about general graphs.

\begin{conj} Suppose that the edges of a $2d$-regular graph are colored such that the edges of each color form a cycle of length $2d$. Then there is a perfect matching containing one edge of each color.
\end{conj}

For $d=2$, there is a simple proof by Zolt\'an Kir\'aly \cite{K}, who also invented the above formulation of the problem. We do not have strong evidence for this conjecture to be true, but it is mathoverflow-hard.

%\section*{Remarks and acknowledgment}

\end{document}